\documentclass[12pt]{amsart} 
\usepackage{amsmath, amssymb,graphics}
\usepackage[colorlinks=true]{hyperref}
\usepackage{fullpage}

\newtheorem{hilfssatz}{Hilfssatz}
\newtheorem{satz}{Satz}
 
\newtheorem{lemma}{Lemma}
\newtheorem{theorem}{Theorem}

\theoremstyle{remark}
\newtheorem*{note-nonum}{Note}

\begin{document}

\title{Ein Satz \"{U}ber die Diophantische Gleichung $Ax^{2}-By^{4}=C$ ($C=1,2,4)$}
\footnote{typeset into \LaTeX\, by N. Tzanakis and P, Voutier}

\author{Wilhelm Ljunggren}

\maketitle

\tableofcontents

\section{Introduction}

In this document, we provide both the original German version of Ljunggren's
article, ``Ein Satz \"{u}ber die Diophantische Gleichung $Ax^{2}-By^{4}=C$ ($C=1,2,4$)'',
Tolfte Skand.~Matematikerkongressen, Lund, 1953, pp.~188--194 (1954), as well
as our English translation.

In the English translation, some footnotes clarifying a few things in the proof
of his Hilfssatz~\ref{lem:Ger-2} are also provided, along with a correction to the statement
of his Satz~\ref{thm:Ger-4}. A justification of this correction is provided in Section~4.
Our next revision of this document will also contain a more detailed explanation
of his proofs.

As this paper of Ljunggren is
difficult to obtain, even in its native German, we hope that this is also of
interest and benefit to researchers in this area.

\newpage

\section{Original German text}

Es seien $A$, $B$ und $C$ ganz-rationale, positive Zahlen und $C=1$, $2$ oder
$4$. Ist $C$ gerade, dann soll $AB$ ungerade sein. Es sei weiter $A$ quadratfrei,
$AB$ keine Quadratzahl und $C=2$ f\"{u}r $A=1$. In einer fr\"{u}heren Arbeit
\cite{Ger-L1} habe ich den folgenden Satz bewiesen:

\begin{satz}
\label{thm:Ger-I}
Die Gleichung
\begin{equation}
\label{eq:Ger-1}
Ax^{2}-By^{4}=C
\end{equation}
hat im Falle $C=1$ h\"{o}chstens eine L\"{o}sung in positiven ganz-rationalen
Zahlen $x$, $y$. Ist $B \equiv -1 \pmod{4}$, so gilt dies auch f\"{u}r $C=4$. In
allen anderen F\"{a}llen gibt es h\"{o}chstens zwei L\"{o}sungen.
\end{satz}

In der citierten Arbeit wurde auch eine Methode angegeben, um die m\"{o}glichen
L\"{o}sungen wirklich zu bestimmen. Diese erfordert die Berechnung der Grundeinheiten
in zwei Zahlk\"{o}rpern vierten Grades. Der Zweck dieser
kleinen Abhandlung ist es eine einfachere L\"{o}sungsmethode nebst eine
Versch\"{a}rfung des obenerw\"{a}hnten Satzes anzugeben.

Man braucht nur solche Werte von $A$, $B$ und $C$ zu betrachten, f\"{u}r welche
eine Gleichung der Form
\begin{equation}
\label{eq:Ger-2}
Az_{1}^{2}-Bz_{2}^{2}=C
\end{equation}
in nat\"{u}rlichen Zahlen $z_{1}$ und $z_{2}$ l\"{o}sbar ist. Dies kann immer
einfach entschieden werden \cite{Ger-A2}.

Zuerst werden wir zwei Hilfss\"{a}tze beweisen, die auch f\"{u}r andere Fragen
in der Theorie der diophantischen Gleichungen vierten Grades von Bedeutung sind.

Es sei $\varepsilon>1$ eine Einheit mit der Norm $+1$ im K\"{o}rper $K \left( \sqrt{D} \right)$,
$D>0$, und es sei weiter $\varepsilon'$ die konjugierte Einheit, also $\varepsilon\varepsilon'=+1$.
Sei $n$ eine ungerade nat\"{u}rliche Zahl. Wir f\"{u}hren folgende Bezeichnungen
ein:
\begin{center}
End of page 188
\end{center}
\hrulefill
\newpage

\begin{align*}
P_{n}(\varepsilon) &= {\varepsilon'}^{\frac{n-1}{2}}\frac{\varepsilon^{n}-1}{\varepsilon-1}
=\frac{\varepsilon^{\frac{n+1}{2}}-{\varepsilon'}^{\frac{n+1}{2}}}{\varepsilon-\varepsilon'}
+\frac{\varepsilon^{\frac{n-1}{2}}-{\varepsilon'}^{\frac{n-1}{2}}}{\varepsilon-\varepsilon'}, \\
Q_{n}(\varepsilon) &= {\varepsilon'}^{\frac{n-1}{2}}\frac{\varepsilon^{n}+1}{\varepsilon+1}
=\frac{\varepsilon^{\frac{n+1}{2}}-{\varepsilon'}^{\frac{n+1}{2}}}{\varepsilon-\varepsilon'}
-\frac{\varepsilon^{\frac{n-1}{2}}-{\varepsilon'}^{\frac{n-1}{2}}}{\varepsilon-\varepsilon'}.
\end{align*}

Man best\"{a}tigt leicht die folgende Relation

\begin{equation}
\label{eq:Ger-3}
P_{n}(\varepsilon) \cdot Q_{n}(\varepsilon)
=\frac{\varepsilon^{n}-{\varepsilon'}^{n}}{\varepsilon-\varepsilon'}.
\end{equation}

\begin{hilfssatz}
\label{lem:Ger-1}
$P_{n}(\varepsilon)$ ist keine Quadratzahl f\"{u}r $n>3$.
\end{hilfssatz}

{\sc Beweis:} Wir m\"{u}ssen zwei Falle unterscheiden:

\vspace*{3.0mm}

{\sc Erster F\"{a}ll:} $n=4t+1$.

Die Gleichung
\[
{\varepsilon'}^{\frac{n-1}{2}} \frac{\varepsilon^{n}-1}{\varepsilon-1}
=z^{2}
\]
kann so geschrieben werden
\[
\left( \varepsilon^{2t} \right)^{2}\varepsilon-z^{2}\varepsilon^{2t} \left( \varepsilon-1 \right)
=1
\]

Hieraus folgt dass die Zahl
\[
\left( \varepsilon^{2t} \sqrt{\varepsilon} + z\varepsilon^{t}\sqrt{\varepsilon-1} \right)^{2}
=-1+2\varepsilon^{4t+1}+2z\varepsilon^{3t}\theta,
\hspace*{3.0mm}
\theta=\sqrt{\varepsilon(\varepsilon-1)},
\]
eine Einheit im Ringe $R$ $(1,\varepsilon,\theta, \theta \varepsilon)$ mit der Relativnorm
$+1$ ist. Nach {\it Dirichlet} haben wir hier in $R$ zwei Grundeinheiten, und die
Einheiten mit der Relativnorm $+1$ sind durch die Potenzen von $\left( \sqrt{\varepsilon}
+\sqrt{\varepsilon+1} \right)$ gegeben. Vgl \cite{Ger-L3} S.8, 12. Dann haben wir
\[
\varepsilon^{2t}\sqrt{\varepsilon} + z\varepsilon^{t}\sqrt{\varepsilon-1}
= \left( \sqrt{\varepsilon}+\sqrt{\varepsilon+1} \right)^{m} = \lambda^{m},
\hspace*{3.0mm}
\text{($m>0$ und ungerade). Mit}
\]
den Bezeichnungen
\[
\lambda' = \sqrt{\varepsilon}-\sqrt{\varepsilon-1}, \hspace*{3.0mm}
\lambda'' = \sqrt{\varepsilon'}+\sqrt{\varepsilon'+1}, \hspace*{3.0mm}
\text{ und } \hspace*{3.0mm}
{\lambda''}' = \sqrt{\varepsilon'}-\sqrt{\varepsilon'-1},
\]
finden wir leicht
\begin{equation}
\label{eq:Ger-4}
\left( \lambda^{m}+{\lambda'}^{m} \right) \left( {\lambda''}^{m}+{{\lambda''}'}^{m} \right)=4,
\hspace*{3.0mm}
\left( \lambda+\lambda' \right) \left( \lambda''+{\lambda''}' \right)=4.
\end{equation}

Hieraus schlie{\ss}en wir
\[
\left( \lambda \lambda'' \right)^{m} + \left( \lambda' {\lambda''}' \right)^{m} -2
=2x^{2}i\sqrt{N}
\]
\[
\lambda \lambda'' + \lambda' {\lambda''}' - 2 = 2i\sqrt{N},
\vspace*{5.0mm}
N=\varepsilon+\varepsilon'-2>0.
\]

\begin{center}
End of page 189
\end{center}
\hrulefill
\newpage

Aus der Formel
\[
\frac{p^{m}+q^{m}-2}{p+q-2}
= \sum_{h=0}^{h=m-1} \frac{m}{h+1} \binom{m+h}{2h+1}(p+q-2)^{h},
\hspace*{3.0mm} pq=1,
\]
\cite{Ger-L3}, S. 15, erh\"{a}lt man f\"{u}r $p=\lambda\lambda''$, $q=\lambda' {\lambda''}'$
\begin{align}
\label{eq:Ger-5}
x^{2} &= m^{2}+\frac{m^{2}(m^{2}-1^{2})}{4!}2\left( 2i\sqrt{N} \right)
+\frac{m^{2}(m^{2}-1^{2})(m^{2}-2^{2})}{6!}2\left( 2i\sqrt{N} \right)^{2}
+ \cdots + \\
&+ \frac{m^{2}(m^{2}-1^{2})(m^{2}-2^{2}) \cdots (m^{2}-k^{2})}{(2k+2)!} 2\left( 2i\sqrt{N} \right)^{k}
+\cdots + \left( 2i\sqrt{N} \right)^{m-1}. \nonumber
\end{align}

Weil $x^{2}$ reell ist, muss der imagin\"{a}re Teil der rechten Seite dieser
Gleichung verschwinden. Nach K\"{u}rzung mit $4\sqrt{N}$ findet sich dann
\[
0=\frac{m^{2}(m^{2}-1^{2})}{4!} - \frac{m^{2}(m^{2}-1^{2})(m^{2}-2^{2})(m^{2}-3^{2})}{8!} 4N + \cdots.
\]

In \cite{Ger-L3}, S. 20 ist gezeigt worden, dass eine solche Gleichung nur f\"{u}r
$m=1$, ($m>0$), erf\"{u}llt ist. Also ist $t=0$ und $n=1$.

\vspace*{3.0mm}

{\sc Zweiter F\"{a}ll:} $n=4t+3$.

Die oben erw\"{a}hnte Methode kann hier nicht angewandt werden um eine obere
Schranke f\"{u}r $n$ zu bestimmen, weil wir in diesem Falle im Ringe $R$
($1$, $\varepsilon$, $\theta_{1}$, $\varepsilon\theta_{1}$), $\theta_{1}=\sqrt{\varepsilon-1}$,
arbeiten m\"{u}ssen, wo die Grundeinheit mit der Relativnorm $+1$ nicht bekannt
ist. Doch k\"{o}nnen wir auch in diesem Falle beweisen, dass es h\"{o}chstens eine
L\"{o}sung gibt. Wir bekommen n\"{a}mlich auch hier ein System der Form~\eqref{eq:Ger-4},
wo aber die zweite Gleichung aus der ersten hergeleitet werden kann.

Aus der Gleichung $P_{n}(\varepsilon)=z^{2}$ ergibt sich
\begin{align}
\label{eq:Ger-6}
\frac{\varepsilon^{2t+2}-{\varepsilon'}^{2t+2}}{\varepsilon-\varepsilon'}
+ \frac{\varepsilon^{t+1}-{\varepsilon'}^{t+1}}{\varepsilon-\varepsilon'}
\left( \varepsilon^{t}+{\varepsilon'}^{t} \right)
&= z^{2}+1, \text{ oder}\\
\frac{\varepsilon^{t+1}-{\varepsilon'}^{t+1}}{\varepsilon-\varepsilon'}
\left( \varepsilon^{t+1}+{\varepsilon'}^{t+1}+\varepsilon^{t}+{\varepsilon'}^{t} \right)
&= z^{2}+1. \nonumber
\end{align}

Zuerst zeigen wir, dass wir $t \equiv 0 \pmod{3}$ haben m\"{u}ssen. Im Falle
$t \equiv 2 \pmod{3}$ ist der erste Faktor der linken Seite von \eqref{eq:Ger-6}
teilbar durch $\dfrac{\varepsilon^{3}-{\varepsilon'}^{3}}{\varepsilon-{\varepsilon'}}
=\left( \varepsilon+\varepsilon' \right)^{2}-1 \equiv 0$ oder $-1 \pmod{4}$. Dies
ist aber unm\"{o}glich, weil die rechte Seite nur ungerade Primfaktoren der
Form $4k+1$ enth\"{a}lt und durch $4$ nicht teilbar ist. Im F\"{a}lle $t \equiv 1 \pmod{3}$
verschwindet die Klammer f\"{u}r $\varepsilon=\varepsilon'=-1$ und f\"{u}r
$\varepsilon=-\varrho$, $\varepsilon'=-\varrho^{2}$, wo $\varrho$ eine primitive
dritte Einheitwurzel

\begin{center}
End of page 190
\end{center}
\hrulefill
\newpage

dedeutet. Die Klammer ist folglich teilbar durch $\varepsilon+\varepsilon'+2$
und durch $\varepsilon+\varepsilon'-1$. Hieraus ergibt sich $\varepsilon+\varepsilon'
\equiv 3 \pmod{8}$ als einzige M\"{o}glichkeit. Die Zahl $t$ muss gerade sein;
sonst tritt n\"{a}mlich auch $\varepsilon+\varepsilon'$ als Teiler auf. Der erste
Faktor kann dann so geschrieben werden, $t=2h$,
\[
\frac{\varepsilon^{t+1}-{\varepsilon'}^{t+1}}{\varepsilon-\varepsilon'}
= -1 + \frac{\varepsilon^{h+1}-{\varepsilon'}^{h+1}}{\varepsilon-\varepsilon'}
\left( \varepsilon^{h}-{\varepsilon'}^{h} \right)
\equiv -1 \pmod{4},
\]
weil $\dfrac{\varepsilon^{h+1}-{\varepsilon'}^{h+1}}{\varepsilon-\varepsilon'}$
durch $\dfrac{\varepsilon^{3}-{\varepsilon'}^{3}}{\varepsilon-\varepsilon'}
=\left( \varepsilon+\varepsilon' \right)^{2}-1$, $\equiv 0 \pmod{4}$, teilbar ist,
Der F\"{a}ll $t \equiv 1 \pmod{3}$ ist folglich auch ausgeschlossen.

Wir k\"{o}nnen nun setzen $n=3^{r}m$, $(m,3)=1$. Die Gleichung $P_{n}(\varepsilon)=z^{2}$
gibt dann
\begin{equation}
\label{eq:Ger-7}
P_{m}(\varepsilon) P_{3^{r}} \left( \varepsilon^{m} \right)=z^{2}.
\end{equation}

Der gr\"{o}sste gemeinsame Teiler der beiden Faktoren ist ein Teiler von $3^{r}$,
$P_{m}(\varepsilon)$ ist aber nicht durch $3$ teilbar. In \cite{Ger-L3}, S.58 ist
n\"{a}mlich gezeigt worden, dass dies f\"{u}r $\dfrac{\varepsilon^{m}-{\varepsilon'}^{m}}
{\varepsilon-\varepsilon'}$ und $\dfrac{\left( \varepsilon^{m} \right)^{3^{r}}-\left( {\varepsilon'}^{m} \right)^{3^{r}}}
{\varepsilon^{m}-{\varepsilon'}^{m}}$ gilt, und unsere Behauptung folgt dann aus
\eqref{eq:Ger-3}. Aus \eqref{eq:Ger-7} ergibt sich folglich
\[
P_{m}(\varepsilon)=w_{1}^{2}, \hspace*{3.0mm} P_{3^{r}}\left( \varepsilon^{m} \right)=w^{2}.
\]

Ist $r$ gerade, so muss $m \equiv 3 \pmod{4}$ sein. Die Gleichung $P_{m}(\varepsilon)=w_{1}^{2}$
ist aber unm\"{o}glich f\"{u}r $m \not\equiv 0 \pmod{3}$. Folglich ist $r$ ungerade
und $m \equiv 1 \pmod{4}$. Dies gibt $m=1$ und
\begin{equation}
\label{eq:Ger-8}
P_{3^{r}}(\varepsilon)=w^{2}.
\end{equation}

Aus \eqref{eq:Ger-8} folgt
\[
\prod_{i=0}^{r-1} P_{3} \left( \varepsilon^{3^{i}} \right)=w^{2}.
\]

Wie in \cite{Ger-L3} S.59 schlie{\ss}t man so
\[
P_{3} \left( \varepsilon^{3^{i}} \right)=w_{i}^{2},
\hspace*{3.0mm} i=0,1,2,3,\ldots, r-1.
\]

Im F\"{a}lle $r>1$ ist dann $P_{3} \left( \varepsilon^{3} \right)=w_{1}^{2}$ oder
$\varepsilon^{3}+{\varepsilon'}^{3}+1=w_{1}^{2}$. Diese letzte Gleichung kann
auch in der folgenden Weise geschrieben werden
\begin{equation}
\label{eq:Ger-9}
\left( \varepsilon+\varepsilon'+2 \right)
\left( \varepsilon+\varepsilon'-1 \right)^{2}
=w_{1}^{2}+1.
\end{equation}

Dies ist aber unm\"{o}glich, weil die rechte Seite entweder $\equiv 0 \pmod{4}$
ist oder einen Primfaktor der Form $4k-1$ enth\"{a}lt. Also ist $r=1$ und $n=3$.
Unser Hilfssatz~\ref{lem:Ger-1} ist damit bewiesen.

\begin{center}
End of page 191
\end{center}
\hrulefill
\newpage

\begin{hilfssatz}
\label{lem:Ger-2}
$Q_{n}(\varepsilon)$ ist keine Quadratzahl f\"{u}r $n>3$ falls $\varepsilon+\varepsilon'$
ungerade ist.
\end{hilfssatz}

{\sc Beweis:} Wir unterscheiden auch hier zwei F\"{a}lle.

{\sc Erster F\"{a}ll:} $n=4t+1$. Die Gleichung $Q_{n}(\varepsilon)=z^{2}$ gibt
\begin{equation}
\label{eq:Ger-10}
\left( \varepsilon^{t}+{\varepsilon'}^{t} \right)
\left( \frac{\varepsilon^{t+1}-{\varepsilon'}^{t+1}}{\varepsilon-\varepsilon'}
- \frac{\varepsilon^{t}-{\varepsilon'}^{t}}{\varepsilon-\varepsilon'} \right)
=z^{2}+1.
\end{equation}

Diese Gleichung ist erf\"{u}llt f\"{u}r $t=0$ ($n=1$). Ist $t>0$ und gerade, so
setzen wir $t=2^{p}t_{1}$, $\left( 2,t_{1} \right)=1$, $p\geq 1$. Dann ist
$\varepsilon^{2^{p}}+{\varepsilon'}^{2^{p}} \equiv -2 \pmod{8}$ ein Teiler von
$\varepsilon^{t}+{\varepsilon'}^{t}$, was aber unm\"{o}glich ist. Also muss $t$
ungerade sein. Im F\"{a}lle $t \equiv 0 \pmod{3}$ ist dann $\dfrac{\varepsilon^{3}+{\varepsilon'}^{3}}{\varepsilon+\varepsilon'}
=\left( \varepsilon+\varepsilon' \right)^{2}-3 \equiv -2 \pmod{8}$ ein Teiler
von $z^{2}+1$. Dies ist aber auch unm\"{o}glich. Ist $t \equiv 1 \pmod{3}$ so
sind $\varepsilon+\varepsilon'$ und $\varepsilon+\varepsilon'-1$ gleichzeitig
Teiler von $z^{2}+1$, was wieder unm\"{o}glich ist. Also muss $t \equiv 2 \pmod{3}$,
d.h. $n \equiv 0 \pmod{3}$.

\vspace*{3.0mm}

{\sc Zweiter F\"{a}ll:} $n=4t+3$. Die Gleichung $Q_{n}(\varepsilon)=z^{2}$ gibt nun
\begin{equation}
\label{eq:Ger-11}
\left( \varepsilon^{t+1}+{\varepsilon'}^{t+1} \right)
\left( \frac{\varepsilon^{t+1}-{\varepsilon'}^{t+1}}{\varepsilon-\varepsilon'}
- \frac{\varepsilon^{t}-{\varepsilon'}^{t}}{\varepsilon-\varepsilon'} \right)
=z^{2}+1.
\end{equation}

Wie im ersten F\"{a}lle ergibt sich auch hier $n \equiv 0 \pmod{3}$.

In beiden F\"{a}llen\footnote{PV correction: ``F\"{a}lle'' and ``F\"{a}llen'' are
used throughout. Sometimes in the collected works, the first one appeared incorrectly
as ``Falle''.} haben wir somit $n=3^{r}m$, $(m,3)=1$, $r \geq 1$ f\"{u}r
$n>1$. Ist $r \geq 1$ und ungerade, so bekommt man, genau wie fr\"{u}her, $n=3$.
Ist $r>1$ und gerade, so muss man auch die Gleichung $\varepsilon^{3}+{\varepsilon'}^{3}-1
=3w_{1}^{2}$ untersuchen. Schreiben wir aber diese Gleichung in der Form
\[
\left( \varepsilon+{\varepsilon'} \right)^{3}-3\left( \varepsilon+\varepsilon' \right)-1
=3w_{1}^{2},
\]
so sehen wir leicht, dass sie mod $9$ unm\"{o}glich ist. Unser Hilfssatz~\ref{lem:Ger-2}
ist damit bewiesen.

\vspace*{5.0mm}

Nun betrachten wir unsere Gleichung~\eqref{eq:Ger-1}. Sind $z_{1}=a$, $z_{2}=b$
die kleinsten positiven L\"{o}sungen von \eqref{eq:Ger-2}, und setzen wir ferner
\[
\varepsilon=\left( \frac{a\sqrt{A}+b\sqrt{B}}{\sqrt{C}} \right)^{2}
=\frac{2Bb^{2}+C+2ab\sqrt{AB}}{C},
\]
so haben wir
\begin{equation}
\label{eq:Ger-12}
\frac{x\sqrt{A}+y^{2}\sqrt{B}}{\sqrt{C}}
=\left( \frac{a\sqrt{A}+b\sqrt{B}}{\sqrt{C}} \right)^{n},
\hspace*{3.0mm} \text{$n>0$ und ungerade, d.h.}
\end{equation}

\begin{equation}
\label{eq:Ger-13}
y^{2}=b\frac{\varepsilon^{\frac{n}{2}}-{\varepsilon'}^{\frac{n}{2}}}{\varepsilon^{\frac{1}{2}}-{\varepsilon'}^{\frac{1}{2}}}
=bP_{n}(\varepsilon).
\end{equation}

\begin{center}
End of page 192
\end{center}
\hrulefill
\newpage

Wir k\"{o}nnen nun $b=rk^{2}$ setzen, wo $r$ ohne quadratische Factoren ist. Aus
\eqref{eq:Ger-13}
folgt
\begin{equation}
\label{eq:Ger-14}
P_{n}(\varepsilon)=rh^{2}.
\end{equation}

Es ist weiter leicht einzusehen, dass $r$ ein Teiler von $n$ ist, d.h. $n=rn_{1}$
und $r$ ungerade. Der Gleichung \eqref{eq:Ger-14} k\"{o}nnen wir dann die folgende
Gestalt geben
\begin{equation}
\label{eq:Ger-15}
P_{n_{1}}(\varepsilon)P_{r} \left( \varepsilon^{n_{1}} \right)=rh^{2}.
\end{equation}

Hieraus ergibt sich
\begin{equation}
\label{eq:Ger-16}
P_{n_{1}}(\varepsilon)=h_{1}^{2} \hspace*{3.0mm}
\text{ und } \hspace*{3.0mm}
P_{r} \left( \varepsilon^{n_{1}} \right)=rh_{2}^{2}.
\end{equation}

Der gr\"{o}sste gemeinsame Teiler der beiden Faktoren der linken Seite von
\eqref{eq:Ger-15} geht n\"{a}mlich in $r$ auf, und ferner ist der zweite dieser
Faktoren $\equiv r \pmod{r^{2}}$.

Aus \eqref{eq:Ger-16} ergibt sich $n_{1}=1$ oder $n_{1}=3$, dem Hilfssatze~\ref{lem:Ger-1}
zufolge. F\"{u}r den Exponenten $n$ in \eqref{eq:Ger-12} haben wir also h\"{o}chstens
zwei M\"{o}glichkeiten:
$n=r$ oder $=3r$. {\it Die Gleichung~\eqref{eq:Ger-1} hat folglich h\"{o}chstens
zwei ganzzahlige L\"{o}sungen $x$, $y$.} Aus \eqref{eq:Ger-16} folgt {\it die notwendige
Bedingung $\varepsilon+\varepsilon'+1=h_{1}^{2}=\dfrac{4}{C}Bb^{2}+3$ f\"{u}r
die Existenz von zwei L\"{o}sungen.} Diese Bedingung ist nie erf\"{u}llt f\"{u}r
$C=1$ oder f\"{u}r $C=4$ mit $B \equiv -1 \pmod{4}$ oder f\"{u}r $C=2$ mit
$B \equiv 1 \pmod{4}$. Damit haben wir Satz~\ref{thm:Ger-I} aufs neue bewiesen, und
gleichzeitig ist eine einfachere Methode angegeben, um die m\"{o}glichen L\"{o}sungen
wirklich zu bestimmen. Ausserdem folgt

\begin{satz}
\label{thm:Ger-2}
Ist die Gr\"{o}sse $\dfrac{4}{C}Bb^{2}+3$ keine Quadratzahl, dann hat $Ax^{2}-By^{4}=C$
h\"{o}chstens eine L\"{o}sung in ganzen positiven Zahlen $x$ und $y$.
\end{satz}

\begin{satz}
\label{thm:Ger-3}
Die Gleichung $Ax^{2}-By^{4}=4$ hat h\"{o}chstens eine L\"{o}sung in ganzen positiven
und teilerfremden Zahlen $x$ und $y$.
\end{satz}

Auf \"{a}hnliche Weise finden wir, unter Benutzung des zweiten Hilfssatzes
\footnote{PV correction: was Hilssatzes}:

\begin{satz}
\label{thm:Ger-4}
Die Gleichung $Ax^{4}-By^{2}=4$ hat h\"{o}chstens zwei L\"{o}sungen in ganzen positiven
Zahlen $x$ und $y$. Ist die Gr\"{o}sse $Bb^{2}+1$ keine Quadratzahl, dann gibt es
h\"{o}chstens eine solche L\"{o}sung. Dies gilt auch f\"{u}r L\"{o}sungen in
teilerfremden Zahlen $x$ und $y$.
\end{satz}

\begin{center}
End of page 193
\end{center}
\hrulefill

\newpage

\setcounter{equation}{0}

\section{English translation}
\footnote{translated by N. Tzanakis and P. Voutier}

Let $A$, $B$ and $C$ be rational, positive integers and $C=1$, $2$ or $4$. If $C$
is even, then $AB$ should be odd. Further, let $A$ be square-free, $AB$ not a
square and $C= 2$ for $A=1$. In an earlier work \cite{Eng-L1}, I have proven the following theorem:

\begin{theorem}
\label{thm:Eng-I}
The equation
\begin{equation}
\label{eq:Eng-1}
Ax^{2}-By^{4}=C
\end{equation}
in the case of $C=1$ has at most one solution in positive, rational integers
$x$, $y$. If $B \equiv -1 \pmod{4}$, this also applies to $C=4$. In all other cases
there are at most two solutions.
\end{theorem}

In the work cited, a method was also given to actually determine the possible
solutions. This requires the calculation of the basic units in two number fields
of the fourth degree. The purpose of this little treatise is to provide a simpler
solution method, along with a tightening of the above-mentioned theorem.

One only needs to consider values of $A$, $B$ and $C$ for which an equation of
the form
\begin{equation}
\label{eq:Eng-2}
Az_{1}^{2}-Bz_{2}^{2}=C
\end{equation}
can be solved in natural numbers $z_{1}$ and $z_{2}$. This can always be easily
decided \cite{Eng-A2}.

First we shall prove two lemmas that are also relevant to other questions in
the theory of Diophantine equations of the fourth degree.

Let $\varepsilon>1$ be a unit with norm $+1$ in the field $K \left( \sqrt{D} \right)$,
$D>0$, and further let $\varepsilon'$ be the conjugate unit, so $\varepsilon \varepsilon' =+1$.
Let $n$ be an odd natural number. We introduce the following terms:
\begin{center}
End of page 188
\end{center}
\hrulefill
\newpage

\begin{align*}
P_{n}(\varepsilon) &= {\varepsilon'}^{\frac{n-1}{2}}\frac{\varepsilon^{n}-1}{\varepsilon-1}
=\frac{\varepsilon^{\frac{n+1}{2}}-{\varepsilon'}^{\frac{n+1}{2}}}{\varepsilon-\varepsilon'}
+\frac{\varepsilon^{\frac{n-1}{2}}-{\varepsilon'}^{\frac{n-1}{2}}}{\varepsilon-\varepsilon'}, \\
Q_{n}(\varepsilon) &= {\varepsilon'}^{\frac{n-1}{2}}\frac{\varepsilon^{n}+1}{\varepsilon+1}
=\frac{\varepsilon^{\frac{n+1}{2}}-{\varepsilon'}^{\frac{n+1}{2}}}{\varepsilon-\varepsilon'}
-\frac{\varepsilon^{\frac{n-1}{2}}-{\varepsilon'}^{\frac{n-1}{2}}}{\varepsilon-\varepsilon'}.
\end{align*}

It is easy to confirm the following relation
\begin{equation}
\label{eq:Eng-3}
P_{n}(\varepsilon) \cdot Q_{n}(\varepsilon)
=\frac{\varepsilon^{n}-{\varepsilon'}^{n}}{\varepsilon-\varepsilon'}.
\end{equation}

\begin{lemma}
\label{lem:Eng-1}
$P_{n}(\varepsilon)$ is not a square for $n>3$.
\end{lemma}

{\sc Proof:} We have to distinguish two cases:

\vspace*{3.0mm}

{\sc First Case:} $n=4t+1$.

The equation
\[
{\varepsilon'}^{\frac{n-1}{2}} \frac{\varepsilon^{n}-1}{\varepsilon-1}
=z^{2}
\]
can be written as
\[
\left( \varepsilon^{2t} \right)^{2}\varepsilon-z^{2}\varepsilon^{2t} \left( \varepsilon-1 \right)
=1
\]

From this, it follows that the number
\[
\left( \varepsilon^{2t} \sqrt{\varepsilon} + z\varepsilon^{t}\sqrt{\varepsilon-1} \right)^{2}
=-1+2\varepsilon^{4t+1}+2z\varepsilon^{3t}\theta,
\hspace*{3.0mm}
\theta=\sqrt{\varepsilon(\varepsilon-1)},
\]
is a unit in the ring $R$ $(1,\varepsilon,\theta, \theta \varepsilon)$ with relative
norm $+1$. According to {\it Dirichlet} we have two fundamental units in $R$,
and the units with relative norm $+1$ are given by the powers of $\left( \sqrt{\varepsilon}
+\sqrt{\varepsilon+1} \right)$. Cf. \cite{Eng-L3} p.8, 12. Then we have
\[
\varepsilon^{2t}\sqrt{\varepsilon} + z\varepsilon^{t}\sqrt{\varepsilon-1}
= \left( \sqrt{\varepsilon}+\sqrt{\varepsilon+1} \right)^{m} = \lambda^{m},
\hspace*{3.0mm}
\text{($m>0$ and odd). With}
\]
the notation
\[
\lambda' = \sqrt{\varepsilon}-\sqrt{\varepsilon-1}, \hspace*{3.0mm}
\lambda'' = \sqrt{\varepsilon'}+\sqrt{\varepsilon'+1} \hspace*{3.0mm}
\text{ and } \hspace*{3.0mm}
{\lambda''}' = \sqrt{\varepsilon'}-\sqrt{\varepsilon'-1},
\]
we easily find that
\begin{equation}
\label{eq:Eng-4}
\left( \lambda^{m}+{\lambda'}^{m} \right) \left( {\lambda''}^{m}+{{\lambda''}'}^{m} \right)=4,
\hspace*{3.0mm}
\left( \lambda+\lambda' \right) \left( \lambda''+{\lambda''}' \right)=4.
\end{equation}

We conclude from this that
\[
\left( \lambda \lambda'' \right)^{m} + \left( \lambda' {\lambda''}' \right)^{m} -2
=2x^{2}i\sqrt{N}
\]
\[
\lambda \lambda'' + \lambda' {\lambda''}' - 2 = 2i\sqrt{N},
\vspace*{5.0mm}
N=\varepsilon+\varepsilon'-2>0.
\]

\begin{center}
End of page 189
\end{center}
\hrulefill
\newpage

From the formula
\[
\frac{p^{m}+q^{m}-2}{p+q-2}
= \sum_{h=0}^{h=m-1} \frac{m}{h+1} \binom{m+h}{2h+1}(p+q-2)^{h},
\hspace*{3.0mm} pq=1,
\]
\cite{Eng-L3}, S. 15, is obtained for $p=\lambda\lambda''$, $q=\lambda' {\lambda''}'$
\begin{align}
\label{eq:Eng-5}
x^{2} &= m^{2}+\frac{m^{2}(m^{2}-1^{2})}{4!}2\left( 2i\sqrt{N} \right)
+\frac{m^{2}(m^{2}-1^{2})(m^{2}-2^{2})}{6!}2\left( 2i\sqrt{N} \right)^{2}
+ \cdots + \\
&+ \frac{m^{2}(m^{2}-1^{2})(m^{2}-2^{2}) \cdots (m^{2}-k^{2})}{(2k+2)!} 2\left( 2i\sqrt{N} \right)^{k}
+\cdots + \left( 2i\sqrt{N} \right)^{m-1}. \nonumber
\end{align}

Because $x^{2}$ is real, the imaginary part of the right hand side of this
equation must vanish. After shortening with $4\sqrt{N}$ you will find
\[
0=\frac{m^{2}(m^{2}-1^{2})}{4!} - \frac{m^{2}(m^{2}-1^{2})(m^{2}-2^{2})(m^{2}-3^{2})}{8!} 4N + \cdots.
\]

In \cite{Eng-L3}, S. 20 it was shown that such an equation is only satisfied for
$m=1$, ($m>0$). So $t=0$ and $n=1$.

\vspace*{3.0mm}

{\sc Second Case:} $n=4t+3$.

The above-mentioned method cannot be used here to determine an upper bound for
$n$, because in this case we have to work in the ring $R$ ($1$, $\varepsilon$,
$\theta_{1}$, $\varepsilon\theta_{1}$), $\theta_{1}=\sqrt{\varepsilon-1}$, where
the fundamental unit with the relative norm $+1$ is not known. But even in this
case we can prove that there is at most one solution. Here too we get a system
of the form~\eqref{eq:Eng-4}, but where the second equation can be derived from the
first.

From the equation $P_{n}(\varepsilon)=z^{2}$ results
\begin{align}
\label{eq:Eng-6}
\frac{\varepsilon^{2t+2}-{\varepsilon'}^{2t+2}}{\varepsilon-\varepsilon'}
+ \frac{\varepsilon^{t+1}-{\varepsilon'}^{t+1}}{\varepsilon-\varepsilon'}
\left( \varepsilon^{t}+{\varepsilon'}^{t} \right)
&= z^{2}+1, \text{ or}\\
\frac{\varepsilon^{t+1}-{\varepsilon'}^{t+1}}{\varepsilon-\varepsilon'}
\left( \varepsilon^{t+1}+{\varepsilon'}^{t+1}+\varepsilon^{t}+{\varepsilon'}^{t} \right)
&= z^{2}+1. \nonumber
\end{align}

First we must show that we have $t \equiv 0 \pmod{3}$. In the case of $t \equiv 2 \pmod{3}$,
the first factor on the left-hand side of \eqref{eq:Eng-6} is divisible by
$\dfrac{\varepsilon^{3}-{\varepsilon'}^{3}}{\varepsilon-{\varepsilon'}}
=\left( \varepsilon+\varepsilon' \right)^{2}-1 \equiv 0$ or $-1 \pmod{4}$. But
this is impossible because the right hand side only contains odd prime factors
of the form $4k+1$ and is not divisible by $4$. In the case of $t \equiv 1 \pmod{3}$
the bracket(right word?) $\varepsilon=\varepsilon'=-1$ disappears and for $\varepsilon=-\varrho$,
$\varepsilon'=-\varrho^{2}$, we deduce that $\varrho$ root is a primitive third
root of unity.

\begin{center}
End of page 190
\end{center}
\hrulefill
\newpage

The bracket(right word?) is therefore divisible by $\varepsilon+\varepsilon'+2$ and $\varepsilon+\varepsilon'-1$.
This results in $\varepsilon+\varepsilon' \equiv 3 \pmod{8}$ as the only possibility.
The number $t$ must be even; otherwise $\varepsilon+\varepsilon'$ also appears as
a divisor. The first factor can then be written as $t=2h$,
\[
\frac{\varepsilon^{t+1}-{\varepsilon'}^{t+1}}{\varepsilon-\varepsilon'}
= -1 + \frac{\varepsilon^{h+1}-{\varepsilon'}^{h+1}}{\varepsilon-\varepsilon'}
\left( \varepsilon^{h}-{\varepsilon'}^{h} \right)
\equiv -1 \pmod{4},
\]
because $\dfrac{\varepsilon^{h+1}-{\varepsilon'}^{h+1}}{\varepsilon-\varepsilon'}$
is divisible by $\dfrac{\varepsilon^{3}-{\varepsilon'}^{3}}{\varepsilon-\varepsilon'}
=\left( \varepsilon+\varepsilon' \right)^{2}-1$, $\equiv 0 \pmod{4}$, the case
$t \equiv 1 \pmod{3}$ is consequently also excluded.

We can now set $n=3^{r}m$, $(m,3)=1$. The equation $P_{n}(\varepsilon)=z^{2}$ then gives
\begin{equation}
\label{eq:Eng-7}
P_{m}(\varepsilon) P_{3^{r}} \left( \varepsilon^{m} \right)=z^{2}.
\end{equation}

The greatest common divisor of the two factors is a divisor of $3^{r}$, but $P_{m}(\varepsilon)$
is not divisible by $3$. In \cite{Eng-L3}, S.58 it has been shown that this is true
for $\dfrac{\varepsilon^{m}-{\varepsilon'}^{m}}{\varepsilon-\varepsilon'}$ and
$\dfrac{\left( \varepsilon^{m} \right)^{3^{r}}-\left( {\varepsilon'}^{m} \right)^{3^{r}}}
{\varepsilon^{m}-{\varepsilon'}^{m}}$, and our assertion then follows from \eqref{eq:Eng-3}.
It follows from \eqref{eq:Eng-7} that
\[
P_{m}(\varepsilon)=w_{1}^{2}, \hspace*{3.0mm} P_{3^{r}}\left( \varepsilon^{m} \right)=w^{2}.
\]

If $r$ is even, then $m \equiv 3 \pmod{4}$. However, the equation $P_{m}(\varepsilon)=w_{1}^{2}$
is impossible for $m \not\equiv 0 \pmod{3}$. Hence $r$ is odd and $m \equiv 1 \pmod{4}$.
This gives $m=1$ and
\begin{equation}
\label{eq:Eng-8}
P_{3^{r}}(\varepsilon)=w^{2}.
\end{equation}

From \eqref{eq:Eng-8} follows
\[
\prod_{i=0}^{r-1} P_{3} \left( \varepsilon^{3^{i}} \right)=w^{2}.
\]

As in \cite{Eng-L3} S.59, one concludes that
\[
P_{3} \left( \varepsilon^{3^{i}} \right)=w_{i}^{2},
\hspace*{3.0mm} i=0,1,2,3,\ldots, r-1.
\]

If $r>1$ then $P_{3} \left( \varepsilon^{3} \right)=w_{1}^{2}$ or
$\varepsilon^{3}+{\varepsilon'}^{3}+1=w_{1}^{2}$. This last equation can also
be written in the following way
\begin{equation}
\label{eq:Eng-9}
\left( \varepsilon+\varepsilon'+2 \right)
\left( \varepsilon+\varepsilon'-1 \right)^{2}
=w_{1}^{2}+1.
\end{equation}

But this is impossible because the right-hand side is either $\equiv 0 \pmod{4}$
or contains a prime factor of the form $4k-1$. So $r=1$ and $n=3$. Our Lemma~\ref{lem:Eng-1}
is thus proven.

\begin{center}
End of page 191
\end{center}
\hrulefill
\newpage

\begin{lemma}
\label{lem:Eng-2}
$Q_{n}(\varepsilon)$ is not a square for $n>3$ if $\varepsilon+\varepsilon'$
is odd.
\end{lemma}

{\sc Proof:} Here, too, we distinguish between two cases.

{\sc First Case:} $n=4t+1$. The equation $Q_{n}(\varepsilon)=z^{2}$ gives
\begin{equation}
\label{eq:Eng-10}
\left( \varepsilon^{t}+{\varepsilon'}^{t} \right)
\left( \frac{\varepsilon^{t+1}-{\varepsilon'}^{t+1}}{\varepsilon-\varepsilon'}
- \frac{\varepsilon^{t}-{\varepsilon'}^{t}}{\varepsilon-\varepsilon'} \right)
=z^{2}+1.
\end{equation}

This equation is satisfied for $t=0$ ($n=1$). If $t>0$ and even, we set $t=2^{p}t_{1}$,
$\left( 2,t_{1} \right)=1$, $p\geq 1$. Then $\varepsilon^{2^{p}}+{\varepsilon'}^{2^{p}}
\equiv -2 \pmod{8}$
\footnote{trans.: This is not correct. We get
$\varepsilon^{2^{p}}+{\varepsilon'}^{2^{p}} \equiv 7 \pmod{8}$
when $\varepsilon+\varepsilon'$ is odd. We will prove it inductively.
For $p=1$, we have $\varepsilon^{2}+{\varepsilon'}^{2}=\left( \varepsilon+\varepsilon' \right)^{2}-2\varepsilon\varepsilon'
=\left( \varepsilon+\varepsilon' \right)^{2}-2$. Since $\varepsilon+\varepsilon'$
is odd, it follows that $\left( \varepsilon+\varepsilon' \right)^{2} \equiv 1 \bmod 8$
and so our result holds for $p=1$.

Supposing that the result holds for $p-1$, we consider
$\left( \varepsilon^{2^{p-1}}+{\varepsilon'}^{2^{p-1}} \right)^{2}
=\varepsilon^{2^{p}}+{\varepsilon'}^{2^{p}}+2\varepsilon^{2^{p-1}}{\varepsilon'}^{2^{p-1}}$.
Since $\varepsilon\varepsilon'=1$, we have
$\left( \varepsilon^{2^{p-1}}+{\varepsilon'}^{2^{p-1}} \right)^{2}
=\varepsilon^{2^{p}}+{\varepsilon'}^{2^{p}}+2$.
By our inductive hypothesis,
$\varepsilon^{2^{p-1}}+{\varepsilon'}^{2^{p-1}} \equiv 7 \pmod{8}$, so
$1 \equiv 7^{2} \equiv \varepsilon^{2^{p}}+{\varepsilon'}^{2^{p}}+2 \pmod{8}$.
Hence $\varepsilon^{2^{p}}+{\varepsilon'}^{2^{p}} \equiv 7 \pmod{8}$ follows.}
is a divisor of $\varepsilon^{t}+{\varepsilon'}^{t}$,
\footnote{trans.: not sure if that is true, but the proof that $t$ cannot be even
can be completed as follows:\\
$\varepsilon^{2^{p}}+{\varepsilon'}^{2^{p}}$ is a divisor of $z^{2}+1$
by \eqref{eq:Eng-10}, but that means that $-1$ is a square mod a number congruent
to $3 \pmod{4}$ (since we just saw that
$\varepsilon^{2^{p}}+{\varepsilon'}^{2^{p}} \equiv 7 \pmod{8}$), which is impossible.\\
Note: this fails if $\varepsilon+\varepsilon'$ is even. Then
$\varepsilon^{2^{p}}+{\varepsilon'}^{2^{p}} \equiv 2 \bmod 8$ and $-1$ can be a
square mod such numbers.}
which is impossible. So
$t$ has to be odd. In the case of $t \equiv 0 \pmod{3}$ then $\dfrac{\varepsilon^{3}+{\varepsilon'}^{3}}{\varepsilon+\varepsilon'}
=\left( \varepsilon+\varepsilon' \right)^{2}-3 \equiv -2 \pmod{8}$ is a divisor
of $z^{2}+1$. But this is also impossible. If $t \equiv 1 \pmod{3}$ then $\varepsilon+\varepsilon'$
and $\varepsilon+\varepsilon'-1$ are simultaneously divisors of $z^{2}+1$, which
is again impossible. So $t \equiv 2 \pmod{3}$, i.e. $n \equiv 0 \pmod{3}$.


{\sc Second Case:} $n=4t+3$. The equation $Q_{n}(\varepsilon)=z^{2}$ now gives
\begin{equation}
\label{eq:Eng-11}
\left( \varepsilon^{t+1}+{\varepsilon'}^{t+1} \right)
\left( \frac{\varepsilon^{t+1}-{\varepsilon'}^{t+1}}{\varepsilon-\varepsilon'}
- \frac{\varepsilon^{t}-{\varepsilon'}^{t}}{\varepsilon-\varepsilon'} \right)
=z^{2}+1.
\end{equation}

As in the first case, $n \equiv 0 \pmod{3}$ also results here.

In both cases we have $n=3^{r}m$, $(m,3)=1$, $r \geq 1$ for $n>1$. If $r \geq 1$
and is odd, then one gets, just as before, $n=3$. If $r>1$ and even, then one
must also examine equation $\varepsilon^{3}+{\varepsilon'}^{3}-1=3w_{1}^{2}$.
But if we write this equation in the form
\[
\left( \varepsilon+{\varepsilon'} \right)^{3}-3\left( \varepsilon+\varepsilon' \right)-1
=3w_{1}^{2},
\]
we can easily see that it is impossible mod $9$. Our Lemma~\ref{lem:Eng-2} is thus proven.


Now let's look at our equation~\eqref{eq:Eng-1}. If $z_{1}=a$, $z_{2}=b$ are the smallest
positive solutions of \eqref{eq:Eng-2}, and we further set
\[
\varepsilon=\left( \frac{a\sqrt{A}+b\sqrt{B}}{\sqrt{C}} \right)^{2}
=\frac{2Bb^{2}+C+2ab\sqrt{AB}}{C},
\]
so we have
\begin{equation}
\label{eq:Eng-12}
\frac{x\sqrt{A}+y^{2}\sqrt{B}}{\sqrt{C}}
=\left( \frac{a\sqrt{A}+b\sqrt{B}}{\sqrt{C}} \right)^{n},
\hspace*{3.0mm} \text{$n>0$ and odd, i.e.,}
\end{equation}
\begin{equation}
\label{eq:Eng-13}
y^{2}=b\frac{\varepsilon^{\frac{n}{2}}-{\varepsilon'}^{\frac{n}{2}}}{\varepsilon^{\frac{1}{2}}-{\varepsilon'}^{\frac{1}{2}}}
=bP_{n}(\varepsilon).
\end{equation}

\begin{center}
End of page 192
\end{center}
\hrulefill
\newpage

We can now set $b=rk^{2}$, where $r$ has no square factors. From \eqref{eq:Eng-13} follows
\begin{equation}
\label{eq:Eng-14}
P_{n}(\varepsilon)=rh^{2}.
\end{equation}

It is also easy to see that $r$ is a divisor of $n$, i.e. $n=rn_{1}$ and $r$ is
odd. We can then give equation~\eqref{eq:Eng-14} the following form
\begin{equation}
\label{eq:Eng-15}
P_{n_{1}}(\varepsilon)P_{r} \left( \varepsilon^{n_{1}} \right)=rh^{2}.
\end{equation}

From this it follows that
\begin{equation}
\label{eq:Eng-16}
P_{n_{1}}(\varepsilon)=h_{1}^{2} \hspace*{3.0mm}
\text{ and } \hspace*{3.0mm}
P_{r} \left( \varepsilon^{n_{1}} \right)=rh_{2}^{2}.
\end{equation}

The greatest common divisor of the two factors on the left-hand side of \eqref{eq:Eng-15}
is divided into $r$, and furthermore the second of these factors is $\equiv r \pmod{r^{2}}$.

From \eqref{eq:Eng-16} we get $n_{1}=1$ or $n_{1}=3$, from Lemma~\ref{lem:Eng-1}. For
the exponent $n$ in \eqref{eq:Eng-12} we have at most two options: $n=r$ or $=3r$.
{\it The equation~\eqref{eq:Eng-1} consequently has at most two integer solutions
$x$, $y$.} From \eqref{eq:Eng-16} {\it follows the necessary condition
$\varepsilon + \varepsilon' + 1 = h_{1}^{2} = \dfrac{4}{C}Bb^{2}+3$ for the
existence of two solutions.} This condition is never fulfilled for $C=1$ or for
$C=4$ with $B \equiv -1 \pmod{4}$ or for $C=2$ with $B \equiv 1 \pmod{4}$. With
this we have proven Theorem~\ref{thm:Eng-I} anew, and at the same time a simpler
method is given to actually determine the possible solutions. Also follows

\begin{theorem}
\label{thm:Eng-II}
If the quantity $\dfrac{4}{C}Bb^{2}+3$ is not a square, then $Ax^{2}-By^{4}=C$
has at most one solution in positive integers $x$ and $y$.
\end{theorem}

\begin{theorem}
\label{thm:Eng-III}
The equation $Ax^{2}-By^{4}=4$ has at most one solution in positive and relatively prime
integers $x$ and $y$.
\end{theorem}

We find in a similar way, using the second lemma.

\begin{theorem}
\label{thm:Eng-IV}
The equation $Ax^{4}-By^{2}=4$ has at most two solutions in positive integers
$x$ and $y$. If the quantity $Bb^{2}+1$ is not a square, then there is at most
one such solution. This also applies to solutions in coprime numbers $x$ and $y$.
\footnote{trans.: working through the proof (see the next pages), we need that $b$
in the minimal solution of \eqref{eq:Eng-2} is odd here too.}
\end{theorem}

\begin{center}
End of page 193
\end{center}
\hrulefill

\newpage

\section{Translators' Addendum}

Here we work through the proof of Theorem~\ref{thm:Eng-IV} and demonstrate why we
require the additional assumption that $b$ is odd.

Now let's look at the equation
\begin{equation}
\label{eq:trans-1}
Ax^{4}-By^{2}=C
\tag{trans-1}
\end{equation}

If $z_{1}=a$, $z_{2}=b$ is the smallest solution in positive integers of \eqref{eq:Eng-2},
and we further set
\[
\varepsilon=\left( \frac{a\sqrt{A}+b\sqrt{B}}{\sqrt{C}} \right)^{2}
=\frac{2Bb^{2}+C+2ab\sqrt{AB}}{C},
\]
so we have
\begin{equation}
\label{eq:trans-2}
\frac{x^{2}\sqrt{A}+y\sqrt{B}}{\sqrt{C}}
=\left( \frac{a\sqrt{A}+b\sqrt{B}}{\sqrt{C}} \right)^{n},
\hspace*{3.0mm} \text{$n>0$ and odd.}
\tag{trans-2}
\end{equation}

So we have
\begin{equation}
\label{eq:trans-3}
x^{2}=a\frac{\varepsilon^{\frac{n}{2}}+{\varepsilon'}^{\frac{n}{2}}}
            {\varepsilon^{\frac{1}{2}}+{\varepsilon'}^{\frac{1}{2}}}
=aQ_{n}(\varepsilon),
\tag{trans-3}
\end{equation}
the last equality holding by multiplying the top and bottom of the middle expression
by $\varepsilon^{n/2}$ and using $\varepsilon \varepsilon'=1$.

We can now set $a=rk^{2}$, where $r$ is square-free. From \eqref{eq:trans-3} follows
that
\begin{equation}
\label{eq:trans-4}
Q_{n}(\varepsilon)=rh^{2},
\tag{trans-4}
\end{equation}
for some integer $h$.

It is also easy to see that $r$ is a divisor of $n$, i.e. $n=rn_{1}$ and $r$ is
odd (recall that $n$ is odd, see \eqref{eq:trans-2}). We can then give equation~\eqref{eq:trans-4}
the following form
\begin{equation}
\label{eq:trans-5}
Q_{n_{1}}(\varepsilon)Q_{r} \left( \varepsilon^{n_{1}} \right)=rh^{2}.
\tag{trans-5}
\end{equation}

From this it follows that
\begin{equation}
\label{eq:trans-6}
Q_{n_{1}}(\varepsilon)=h_{1}^{2} \hspace*{3.0mm}
\text{ and } \hspace*{3.0mm}
Q_{r} \left( \varepsilon^{n_{1}} \right)=rh_{2}^{2}.
\tag{trans-6}
\end{equation}

The greatest common divisor of the two factors on the left-hand side of \eqref{eq:trans-5}
is a divisor of $r$, and furthermore the second of these factors is $\equiv r \pmod{r^{2}}$.

From \eqref{eq:trans-6} we get $n_{1}=1$ or $n_{1}=3$, from Lemma~\ref{lem:Eng-2},
provided that $\varepsilon+\varepsilon'$ is odd (i.e., $4Bb^{2}/C+2$ is odd,
so $4Bb^{2}/C$ is odd. This means $C=4$ and that $Bb^{2}$ is odd. We know that
$B$ is odd from our initial assumption that when $C$ is even, like here, then
$AB$ is odd. Thus we need $b$ to be odd.).

For the exponent $n$ in \eqref{eq:trans-2} we have at most two options: $n=r$ or $n=3r$.
{\it The equation~\eqref{eq:trans-2} consequently has at most two integer solutions
$x$, $y$.} From \eqref{eq:trans-6} {\it follows the necessary condition
$\varepsilon + \varepsilon'-1 = h_{1}^{2} = \dfrac{4}{C}Bb^{2}+1$ for the
existence of two solutions.}
\end{document}